\documentclass[12pt]{article}
\usepackage{amsmath,amssymb,amsfonts}
\usepackage{setspace}
  \renewenvironment{thebibliography}[1]{%
    \begin{oldthebibliography}{#1}%
      \setlength{\parskip}{0ex}%
      \setlength{\itemsep}{0ex}%
  }%
  {%
    \end{oldthebibliography}%
  }

\newtheorem{theo}{Theorem}[section]
\newtheorem{pro}[theo]{Proposition}
\newtheorem{lem}[theo]{Lemma}
\newtheorem{cor}[theo]{Corollary}

\newtheorem{rem}[theo]{Remark}

\def\qed{\hbox {\rlap{$\sqcup$}{$\sqcap$}}}

\def\ignore#1{}

\def\bigb{\bigskip\noindent}         
\def\med{\medskip\noindent}
\def\sms{\smallskip}
\def\ms{\medskip}

\def\sm{\smallskip\noindent}
\def\nl{\hfil\break}
\def\ni{\noindent}
\def\noi{\noindent}

\def\today{\noindent\number\day
\space\ifcase\month\or
  January\or February\or March\or April\or May\or June\or
  July\or August\or September\or October\or November\or December\fi
  \space\number\year}

\input amssym.def
\input amssym.tex

\def\bP {{\mathbb P}} \def\bE {{\mathbb E}} 
\def\bR {{\mathbb R}} \def\bN {{\mathbb N}} \def\bZ {{\mathbb Z}}

 \def\sE {{\cal E}} \def\sF {{\cal F}}
  
  \def\sL {{\cal L}}

\def\ol{\overline}

\def\Ra{$\Rightarrow$}

\def\al {\alpha}
\def\lam {\lambda}  \def\Gam{\Gamma}          
   \def\gam{\gamma}                

  \def\eps{\varepsilon}

\def\th{\theta}

\def\pd {\partial}       
\def\q{\quad} \def\qq{\qquad}
\def\dint{\int\kern-.6em\int}

\def \half {{\tfrac12}}

\def \qed {\hfill$\square$\par}

\def\=d{{\,\buildrel (d) \over =\,}}
\def\a.s.{{\buildrel a.s. \over \longrightarrow}}

\def\proof{\sm {\it Proof. }}

\makeatletter
\@addtoreset{equation}{section}

\makeatother

\title{Parabolic Harnack inequality and heat kernel 
estimates for random walks  with long range jumps}

\author{Martin T. Barlow
\thanks{Research partially supported by NSERC (Canada),
the 21st Century COE Program in Kyoto University \indent~ (Japan), 
and by EPSRC (UK).}
\and
Richard F. Bass
\thanks{Research partially supported by
NSF Grant DMS-0601783.}
\and
Takashi Kumagai
\thanks{Research partially supported by the
Grant-in-Aid for Scientific Research (B) 18340027 (Japan).}}

\begin{document}

\maketitle
\date

\begin{abstract}
  We investigate the relationships between the parabolic Harnack
  inequality, heat kernel estimates, some geometric conditions, and
  some analytic conditions for random walks with long range jumps.
  Unlike the case of diffusion processes, the parabolic Harnack
  inequality does not, in general, imply the corresponding heat kernel
  estimates.
\end{abstract}

\section{Introduction}

This paper investigates the relationships between the parabolic
Harnack inequality, heat kernel estimates, some geometric conditions,
and some analytic conditions for random walks with long range jumps.
By random walks with long range jumps, also known as random walks with
unbounded range, we mean random walks for which there does {\bf not}
exist a positive integer $K$ such that the probability of a jump
larger in size than $K$ is zero.

Our investigation combines two lines of research that have received
much attention. For the past few decades there has been a great deal
of interest in extending the results of DeGiorgi, Nash, Moser, and
others on the regularity of solutions to the heat equation on $\bR^d$
with respect to elliptic operators in divergence form to much more
general state spaces.  Among the state spaces considered are
manifolds, graphs, and fractals.  A typical result for the case of
diffusions on manifolds or nearest neighbor random walks on graphs is
along the lines of the following.  (For a precise statement of the
results, see \cite{Gr1,SC1,Del}.)

\begin{theo}\label{ThmA} 
The following are equivalent:
\nl (a) Gaussian upper and lower bounds on the heat kernel;
\nl (b) the parabolic Harnack inequality; 
\nl (c) volume doubling and a family of Poincar\'e inequalities.
\end{theo}

The other line of research leading to this paper is the study of
Harnack inequalities and heat kernel estimates for processes with
jumps on $\bR^d$, $\bZ^d$, or state spaces with similar structure.
These results are more recent; among the early papers here are
\cite{BL, CK}. The motivation is that researchers in mathematical
physics, mathematical finance, and other areas want to allow their
models to have jumps, but many of the basic properties of jump
processes are still unknown. A typical result (see the cited
references for exact statements) is the following.

\begin{theo}\label{ThmB}
If there exists $\alpha\in (0,2)$ such that the probability of a jump
from $x$ to $y$ is comparable to $|x-y|^{-d-\alpha}$, then the
following hold: 
\nl (a) Polynomial type upper and lower bounds on the heat kernel;
\nl (b) the parabolic Harnack inequality.
\end{theo}

It is therefore quite natural to study heat kernel estimates and the
parabolic Harnack inequality for random walks on more general graphs
where there is the possibility of arbitrarily large jumps. Besides
being interesting in its own right, we believe this study sheds
additional light on pure jump processes of all types. It should
also be mentioned that there are significant differences between
the results for the diffusion or nearest neighbor case (Theorem \ref{ThmA})
and the results we obtain here for the long range case (Theorems \ref{Thm1.2}
and \ref{Thm1.3})

In this paper we investigate these connections in the framework
of continuous time random walks on graphs. We believe that our 
results should extend with only minor changes to jump processes on
metric measure spaces. However, in that context 
some issues of regularity would have to be treated.

\ms Let $\Gam=(G,E)$ be an infinite connected graph, where $G$ is the
set of vertices and $E$ the set of edges. We write $d(x,y)$
for the graph distance, and we assume that $\Gam$ is locally finite.
We let $B(x,r)=\{y: d(x,y) \le r\}$ denote balls in the graph metric.
(We allow $r \in [0,\infty)$.)
The notation $x\sim y$ means that $d(x,y)=1$.

Let $J(x,y)=J(y,x)$ be a symmetric non-negative function on $G\times G$
with $J(x,x)=0$ for all $x$.
We write
\begin{equation} J(x,A) =\sum_{y \in A} J(x,y), \label{1.1A}
\end{equation}
and assume there exists $C_J \in [1,\infty)$ such that
\begin{equation}  
C_J^{-1}\le J(x,G) \le C_J, \q x \in G. \label{1.1} 
\end{equation}
Let $\mu$ be a measure on $G$ such that $\mu_x = \mu(\{x\})$
satisfies for some constant $C_M \in [1,\infty)$
\begin{equation} 
C_M^{-1}\le \mu_x \le C_M, \q x \in G.
\label{mucond}
\end{equation} 
We write
\begin{equation} V(x,r) =\mu(B(x,r)). \label{1.2A}
\end{equation} 

We define the quadratic form
\begin{equation}
\sE(f,f) =\half \sum_x \sum_y (f(x)-f(y))^2 J(x,y).  \label{1.2B}
\end{equation}
An easy application of Cauchy-Schwarz shows that 
$\sE(f,f)\le 2C_J C_M ||f||_2^2$. We consider the continuous time Markov
process $X=(X_t, t\ge 0, \bP^x, x\in G)$ with jump rates 
from $x$ to $y$ of $\mu_x^{-1} J(x,y)$.
This is the Markov process associated with the Dirichlet form
$(\sE, L^2(G,\mu))$, and has generator
\begin{equation} 
\sL f(x) = \frac{1}{\mu_x} \sum_y (f(y)-f(x)) J(x,y).   \label{1.2}  
\end{equation}
Note that since $||\sL f||_2 \le 2 C_J C_M ||f||_2$, $\sL$ is defined
on $L^2(G,\mu)$.
We write $p_t(x,y)$ for the heat kernel on $\Gam$; this is the
transition density of the process $X$ with respect to $\mu$:
\begin{equation} 
 p_t(x,y) = \frac{\bP^x( X_t =y)}{\mu_y}. 
\end{equation} 
Since the total jump rate out of $x$ is
$\mu_x^{-1} J(x,G) \le C_J C_M$, the process $X$ is conservative
and $\sum_y p_t(x,y) \mu_y =1 $ for all $x,t$.

\medskip
We now consider various conditions which could be imposed 
on $\Gam$, $J$, $X$ and $p$.

\sm{\it 1. Volume growth.}
$G$ satisfies volume doubling VD if there exists a constant 
$C_V$ such that
$$ V(x,2r) \le C_V V(x,r) \q \hbox { for all } x \in G, 
r>0. \eqno(V\!D)$$
It is easy to check that VD implies that there exists
$\al_1>0$ such that if $x,y \in G$ and $0<r <R$ then
\begin{equation} \label{vd2} 
\frac{V(x,R)}{V(y,r)} 
\le c_1 \Bigl( \frac{d(x,y)+R}{r} \Bigr)^{\al_1}. 
\end{equation}
Thus $V(x,d(x,y)) \asymp V(y,d(x,y))$, where 
$\asymp$ means the ratio of the two sides is bounded above and
below by two positive constants not depending on $x$ or $y$.

VD also implies that
\begin{equation} \label{vd3} 
\frac{V(x,2r)}{V(x,r)} \ge 1 + C_V^{-4}, \q 
\mbox{ provided $r \ge 1$}.
\end{equation}
To see this choose $y \in B(x, 2r)-B(x,r)$, such that
$B(y,r/4) \subset B(x,2r)$ and $B(y,r/4) \cap B(x,r) =\emptyset$.
So, $V(x,2r) \ge V(x,r) + V(y,r/4)$, while by VD, 
$C_V^4 V(y,r/4) \ge V(y,4r) \ge V(x,r)$. Combining these 
gives (\ref{vd3}).

A more restrictive condition is that $V(x,r)$ grows like $r^d$
(where $d\in [1,\infty)$): 
$$ C_1 r^d \le V(x,r) \le C_2 r^d, \q r\ge 1. \eqno(V(d))  $$

\ni{\it 2. Transition density estimates.}
Next we introduce various conditions on the heat kernel $p_t(x,y)$.
$X$ satisfies UHKP$(\al)$ if
$$ p_t(x,y) \leq c_1\frac{1}{V(x,t^{1/\al})} 
\wedge \frac{t}{V(x,R) R^{\al}}, \q t>0, \hbox { where } R=d(x,y),
 \eqno(U\!H\!K\!P(\al)) $$
and LHKP$(\al)$ if
$$ p_t(x,y) \geq c_2\frac{1}{V(x,t^{1/\al})} 
\wedge \frac{t}{V(x,R) R^{\al}}, \q t>0, \hbox { where } R=d(x,y).
 \eqno(L\!H\!K\!P(\al)). $$
If both UHKP$(\al)$ and LHKP$(\al)$ hold, we say HKP$(\al)$ holds.
The `P' here stands for `polynomial' -- this kind of decay arises 
frequently for processes with long range jumps, instead of the 
Gaussian type behavior
associated with continuous processes.  Note that the first term in
UHKP$(\al)$ (and in LHKP$(\al)$) is smaller than the second term 
if and only if $t> R^\al$. 

If we just have the upper bound for $x=y$, we say
UHD$(\al)$ holds:
$$ p_t(x,x) \le \frac{c}{V(x,t^{1/\al})}.  \eqno(U\!H\!D(\al))  $$
If $V(d)$ holds, then  HKP$(\al)$ takes the form:
$$ p_t(x,y) \asymp t^{-d/\al} \wedge \frac{t}{R^{d+\al}}, \q t \ge 0, 
\hbox { where } R=d(x,y). $$

In the proof of the parabolic Harnack inequality a key role is
played by the following lower bound. 
Let $p^{B(x,r)}_t(\cdot,\cdot)$ be the heat kernel of the
process $X$ killed on exiting $B(x,r)$. 
$(\Gam,J)$ satisfies NDLB$(\al)$ if  
there exist constants $c_1, c_2, c_3$ such that 
$$ p^{B(x,r)}_t(x',y') \ge \frac{c_1}{V(x,r)}, 
 \qq x',y' \in B(x,r/2), \q  c_2 r^\al \le t \le c_3 r^\al. 
 \eqno (N\!D\!L\!B(\al)) $$ 

\med {\it 3. Harnack inequalities.}
Let $I$ be an open subset of $\bR$, $A \subset G$, and $Q=I \times A$.
Let $u(t,x)$ be defined on $I \times G$. We say that $u$ is a 
{\sl caloric on $Q$} if
\begin{equation} 
\frac{\pd u}{\pd t}(t,x) = \sL u(t,x), \q t\in I, x\in A. \label{1.6}
\end{equation}
We can interpret (\ref{1.6}) in the weak sense in $t$, that is, for any
$\psi \in C^\infty_c(I)$ 
\begin{equation}  
- \int_I \psi'(t)u(t,x)dt  = \int_I \sL u(t,x)\psi(t) dt. \label{1.7}
\end{equation}
Here $C^\infty_c(I)$ is the set of $C^\infty$ functions with 
compact support
contained in $I$.

\ms 
Let $Q(x_0,T,R)=(0,T)\times B(x_0,R)$,
and set
\begin{align*}
 Q_- &= Q_-(x_0,T,R)=  [T/4, T/2]\times B(x_0, R/2), \\
 Q_+ &=  Q_+(x_0,T,R)= [3 T/4, T ]\times B(x_0, R/2). 
\end{align*}
For $Q$ as above we write $s+Q=(s,s+T)\times B(x_0,R)$.
We say  the parabolic Harnack inequality 
PHI$(\al)$ holds, if for $\lam\in (0,1]$
there exist constants $C_P(\lam)$, depending
only on $\lam$, such that whenever $u=u(t,x)\ge 0$ is 
caloric in $Q(x_0,\lam R^\al,R)$ 
and continuous at time $T$,  then
$$ \sup_{Q_-} u \le C_P(\lam)  \inf_{Q_+} u. \eqno (P\!H\!I(\al))$$
Note that the continuity of $u$ at time $T$ is assumed since we often 
use PHI$(\al)$ at time $T$. Alternatively, one may define PHI$(\al)$ as
``$\cdots$ whenever $u=u(t,x)\ge 0$ is 
caloric in $Q(x_0,\lam R^\al+\eps,R)$ for some $\eps>0$, then $\cdots$ ''. 

We say a function $h$ defined on $G$ is  {\sl harmonic} on a subset $A$ if
$$\sL h(x)=0, \qq x\in A.$$
The elliptic Harnack inequality holds if there exists a constant $c$
not depending on $x_0$ or $R$ 
such that if $h$ is non-negative on $G$ and harmonic in $B(x_0,2R)$,
$x_0\in G$, $R>1$, then 
$$h(x)\leq ch(y), \qq x,y\in B(x_0,R). \eqno (E\!H\!I)$$

\begin{rem}
\begin{itemize}
\setlength{\itemsep}{-\parsep}
{\rm \item[(a)] The classical parabolic Harnack inequality for 
diffusions on manifolds has $\al=2$. Parabolic Harnack inequalities with 
anomalous scaling $\al>2$ are given in {\rm \cite{HSC,BB}}. 
The case $\al<2$ is discussed in {\rm \cite{BL,CK,CK2}}.
\item[(b)] If $R<1$, then $B(x,R)$ and $B(x,\half R)$ are both just 
the single point $\{x\}$. Nevertheless the parabolic Harnack inequality as 
stated above still makes sense, and in fact it is easy to check that this 
local parabolic Harnack inequality,
(i.e., the parabolic Harnack inequality with $R<1$), will
always hold under the condition that $C_J <\infty$.
\item[(c)] If $\al>1$ then the introduction of $\lam$ is 
not necessary, since the parabolic Harnack inequality 
for $\lam=1$ implies the 
parabolic Harnack inequality for any $\lam\in (0,1]$. 
To see why we need to introduce $\lam$ in the case when $\al \le 1$,
let $x_0, x_1 \in G$ with $d(x_0,x_1)=R$, let $0< T < R^\al$ and
suppose we wish to find a chain of  $n$ space-time boxes
$Q_i=s_i+Q(x_i,r^\al,r)$ linking $(x_0,0)$ with  $(x_1,T)$.
If $t=r^{\al}$ then we need $nr \ge R$, and $nt \le T$,
which implies that $n^{\al-1} \ge R^\al/T$. 
Since $n\ge 1$, chaining of this type is only possible when $\al>1$.
\item[(d)] Since a harmonic function is caloric, PHI$(\al)$ implies EHI.}
\end{itemize}
\end{rem}

\ms  
\ni{\it 4. Analytic estimates.} 
Let $P_t^B$ be the semigroup for $X_t$ killed on exiting $B$.
We consider the following semigroup bound: 
there exist $c_1,c_2,c_3>0$ such that 
$$  \|P_t^{B(x_0,r)}\|_{1\to\infty} \le 
c_1 V(x_0,r)^{-1}, 
\q \hbox{ for all } c_2r^\al\le t\le c_3 r^\al, 
r\ge 0, x_0\in G. \eqno(SB(\alpha)) $$ 
For an operator $T$ on functions, $\|T\|_{p\to q}=
\sup\{\|T f\|_q: \|f\|_p\leq 1\}$, where $\|f\|_p$ is the $L^p$ norm 
of $f$ with respect to $\mu$. 
 
\ms
\ni{\it 5. Jump kernel.} We also consider bounds on $J$:
 $$ J(x,y) \le  \frac{c \mu_x \mu_y}{d(x,y)^\al V(x,d(x,y))}, 
\eqno(UJ(\al))$$
 $$J(x,y) \ge  \frac{c \mu_x \mu_y} {d(x,y)^\al V(x,d(x,y))}. 
\eqno(LJ(\al))$$ 
If $J$ satisfies both UJ$(\al)$ and LJ$(\al)$ we say
it satisfies J$(\al)$. 

\nl We will use the following smoothness results on $J$: 
$$ J(x,y) \le \frac{c\mu_x}{V(x,r)} \sum_{y\in B(x,r)} J(x',y) \q
\hbox { whenever $r\le \half d(x,y)$}.  \eqno(U\!J\!S)$$
We say LJS holds if
\begin{align}
\label{ljs1}
 J(x,y) &\ge \frac{c\mu_x}{V(x,r)} \sum_{y\in B(x,r)} J(x',y) \q
\hbox { whenever $r\le \half d(x,y)$}, \\
\label{ljs2}
J(x,y) &\ge c_0>0  \q \mbox{ whenever } x \sim y. 
\end{align} 
(Recall that $x\sim y$ means $d(x,y)=1$.)
We say JS holds if the local non-degeneracy condition (\ref{ljs2})
holds and in addition
\begin{equation}  \label{js}
J(x_1,y) \le c J(x_0,y)\q  \hbox { if } d(x_0,x_1)\le \half d(x_0,y). 
\end{equation} 
Given VD, then combining UJS and LJS gives JS
by a straightforward argument (see Lemma \ref{LemJS}).

\ni {\it 6. Exit times.} For $A \subset G$ we write
$$ \tau_A = \min \{ t\ge 0: X_t \not\in A\}. $$
$(\Gam,J)$ satisfies $E_\al$ if for all 
$x\in G$, $r\ge 1$,
$$ c_1 r^\al \le \bE^x \tau_{B(x,r)} \le c_2 r^\al.  \eqno (E_\al)$$ 
In \cite{GT2} it is proved that for simple random walks
the condition VD+EHI + $E_2$
is also equivalent to conditions (a)--(c) of Theorem \ref{ThmA}.
(The case $\al >2$ is also treated there.)

\ni {\it 7. Poincar\'e inequality.} $(\Gam,J)$ satisfies PI$(\al)$ if 
there exists a constant $C_Q$ such that for any ball $B=B(x,R)\subset G$
with $R\geq 1$ and $f:B\to \bR$,
$$\sum_{x\in B}(f(x)-\ol f_B)^2\mu_x \le 
C_Q R^\al \sum_{x,y\in B} (f(x)-f(y))^2
J(x,y), \eqno (PI(\al))$$
where $\ol f_B=\mu(B)^{-1} \sum_{x\in B} f(x)\mu_x $.

\bigskip The main results of this paper are as follows.
First, we see that some of the implications in Theorem \ref{ThmA}
do hold in the long range jump case. 

\begin{theo}\label{Thm1.1} Let $\al>0$. 
\nl (a)  HKP$(\al)$ implies  PHI$(\al)$.
\nl (b)  PHI$(\al)$ implies VD+ EHI + $E_\al$.
\end{theo}

A counterexample in Section 6  shows that
the converse of  Theorem \ref{Thm1.1}(a) does not hold.
We also sketch an example in that section which shows that
the converse of (b) also fails.
It is easy to see that VD plus PI$(\al)$ is not enough 
to prove PHI$(\al)$.

Given the gap between HKP$(\al)$ and PHI$(\al)$, 
one wishes to find good conditions equivalent to each of these. 

\begin{theo}\label{Thm1.2} 
Assume $V(d)$, and $\al\in (0,2)$. The following are equivalent:
\nl (a)  J$(\al)$,
\nl (b)  HKP$(\al)$,
\nl (c)  PHI$(\al)$ and LJS.
\end{theo}

\begin{theo}\label{Thm1.3} 
Assume $V(d)$ and $\al\in (0,2)$. The following are equivalent:
\nl (a)  PHI$(\al)$
\nl (b)  UJ$(\al)$+ $PI(\al)$ + UJS.
\end{theo}

\begin{rem}\label{Rem-eff} 
{\rm 
1. Theorems \ref{Thm1.2} and 
\ref{Thm1.3} are enough to prove `stability' of
HKP$(\al)$ and PHI$(\al)$ in the following sense. 
We say a property $P$ is stable if whenever 
$J$ and $J'$ are comparable, 
i.e., $J(x,y) \asymp J'(x,y)$ for $x,y\in G$,  
and $P$ holds for $(\Gam,J)$, then $P$ also holds for $(\Gam,J')$.

\noindent 2. Our results are actually slightly stronger than
those stated, in that the constants which arise in the conclusions
only depend on those in the various hypotheses. So, for example,
a more careful statement of Theorem 1.4(a) would be ``Suppose
the graph $\Gam$, jump kernel $J$ and measure $\mu$ satisfy
(\ref{1.1}) and  (\ref{mucond}) with constants $C_J$ and $C_M$,
and that 
$(\Gam,J)$ satisfies  HKP$(\al)$ with constants $C_1$ and $C_2$.
Then $(\Gam,J)$  satisfies  HKP$(\al)$ with a constant 
$C_P$, where $C_P$ depends only on the constants $C_1$, $C_2$,
$C_J$ and $C_M$.''
}
\end{rem}

Section 2 shows that HKP$(\al)$ implies a lower bound on the heat kernel
of the killed process in a ball. Section 3 proves the parabolic Harnack 
inequality, using the `balayage' argument introduced in \cite{BBCK}.
Section 4 looks at consequences of the parabolic Harnack inequality --
see Proposition \ref{phicons} for a summary of these. 
Section 5 looks at consequences of the condition SB$(\al)$,
and combining these with the results of Sections 2--4 
completes the proofs of Theorems  \ref{Thm1.1} -- \ref{Thm1.3}.
In Section 6 we give some counterexamples, which show that
the converses of (a) and (b) in Theorem \ref{Thm1.1} do not hold.

\sm {\bf Note.} By Theorem \ref{Thm1.1} 
each of  PHI$(\al)$ and  HKP$(\al)$
implies VD. Some of the implications in the Theorems \ref{Thm1.2} and
\ref{Thm1.3} do not need $V(d)$; in addition some do not need the
condition $\al <2$.
We summarize this in the following table.

\begin{tabbing}
xxxxxxxxxxxxxxxxxxxxxxxxx\= xxxxxxxxxxxxxxxxxxx\= xxxxxxxxxxxxxxxxxx\=\kill
{\it Statement} \>{\it Volume condition} \> {\it Range of $\al$} \\
~~\>~~ \>~~\\
Theorem \ref{Thm1.2}(a)\Ra (b)   \> $V(d)$   \> $0< \al< 2$ \\
Theorem \ref{Thm1.2}(b)\Ra (c)  \> None   \> $0< \al< \infty$ \\
Theorem \ref{Thm1.2}(c)\Ra (a)     \> None   \> $0< \al< 2$ \\
Theorem \ref{Thm1.3}(a)\Ra (b)  \> None    \>  $0< \al< \infty$ \\
Theorem \ref{Thm1.3}(b)\Ra (a)  \> $V(d)$   \> $0< \al< 2$ 
\end{tabbing}

\ms For the convenience of the reader, we list the abbreviations
we have used and in which subsection of Section 1 they can be found.

\begin{tabbing}
xxxxxxxxxxxxxxxxxxxxxxxxxxxx\= xxxxxxxxxxxxxxxxxxx\= xxxxxxxxxxxxxxx
\= xxxxxxxxxx\=\kill
{\it Abbreviation}\>{\it Meaning} \> \>{\it Subsection}\\
\> \> \> \\
$E_\al$ \> Exit time \> \> 1.6\\
EHI \> Elliptic Harnack inequality \> \> 1.3\\
HKP$(\al)$, LHKP$(\al)$,  UHKP$(\al)$ \> Upper and lower heat kernel \> \> 1.2\\
J$(\al)$, LJ$(\al)$, UJ$(\al)$ \> Jump kernel bounds \> \> 1.5\\
JS, LJS, UJS \> Jump smoothness \> \> 1.5 \\
NDLB$(\al)$ \> Near diagonal lower heat kernel \> \> 1.2\\
PHI$(\al)$ \> parabolic Harnack inequality \> \> 1.3 \\
PI$(\al)$ \>  Poincar${\acute{\rm e}}$ inequality \> \> 1.7 \\
SB$(\al)$ \> 
Semigroup bounds \> \> 1.4 \\ 
UHD$(\al)$ \> On diagonal upper heat kernel \> \> 1.2 \\
VD \> Volume doubling \> \> 1.1 \\
V(d) \> Volume growth \> \> 1.1 
\end{tabbing}

Throughout the paper, we use $c, c'$ to denote strictly positive
finite constants whose values are not significant and may change
from line to line. We write $c_i$ for positive constants whose
values are fixed within theorems and lemmas.
We adopt the convention that if we cite elsewhere the constant $c_1$
in Lemma $2.2$ (for example), we denote it as $c_{2.2.1}$.

\section{The heat kernel killed outside a ball}

We begin by proving that HKP$(\al)$ gives a near-diagonal lower
bound on the heat kernel killed outside a ball.
We assume  HKP$(\al)$ holds with constants $C_1$ and $C_2$, so that,
writing $r=d(x,y)$,  
$$ \frac{C_1}{V(x,t^{1/\al})} \wedge \frac{C_1 t }{V(x,r) r^{\al}} 
\le p_t(x,y) 
\le \frac{C_2}{V(x,t^{1/\al})} \wedge \frac{C_2 t }{V(x,r) r^{\al}}. $$

\begin{lem}\label{Lem2.0} HKP$(\al)$ implies VD. 
\end{lem}

\proof  Let $t^\al =r$. Then
\begin{equation*}
 \frac{C_2}{ V(x, (2t)^{1/\al} ) } \ge p_{2t}(x,x)  
   \ge \sum_{y \in B(x,r) } p_t(x,y)^2 \mu_y 
 \ge V(x,r) \frac { C_1^2 } { V(x,r)^2 } .
\end{equation*}
Rearranging gives $ V(x,2^{1/\al} r) \le c V(x,r)$, 
which implies VD.
\qed

\begin{rem}
{\rm Note that HKP$(\al)$ is equivalent to the following:
\begin{equation}\label{hkp2}
C_1(\frac 1{V(x,C_2t^{1/\al})} \wedge \frac{t}{V(x,r) r^{\al}}) 
\le p_t(x,y) 
\le C_3(\frac 1{V(x,C_4t^{1/\al})} \wedge \frac{t}{V(x,r) r^{\al}}),
\end{equation}
where $r=d(x,y)$. 
Indeed, one can prove VD from (\ref{hkp2}) similarly to the proof of  
Lemma \ref{Lem2.0}, so (\ref{hkp2}) implies HKP$(\al)$.
}
\end{rem}

\begin{lem}\label{Lem2.1} Assume HKP$(\al)$. 
Let $B= B(x_0,R)\subset G$.  For each $\lam\in (0,1)$ there
exists a constant $c_1(\lam)$ such that 
\begin{equation}  
p^B_{\lam R^\al}(x,y) \ge \frac{c_1(\lam)}
{V(x_0,R)},  \q   x,y \in B'=B(x_0,R/2).  
\label{2.1}\end{equation}
\end{lem}
\proof Let $\tau=\tau_B$. 
By the strong Markov property of $X$, 
for $x,y \in B'=B(x_0,R/2)$ and for any $t>0$ 
\begin{equation} 
p_t(x,y) = p^B_t(x,y) + \bE^x 1_{(\tau<t)} p_{t-\tau}(X_\tau,y)
 \le  p^B_t(x,y) + \sup_{0\le s \le t} \sup_{z \in B^c} p_s(z,y). \label{2.2}
\end{equation}
For $z\in B^c$, $y\in B'$, and $s\le (R/2)^\al$ the upper bound
in HKP$(\al)$ gives 
$$ p_s(z,y) \le \frac{C_2 s}{V(x,d(z,y)) d(z,y)^{\al}}
 \le \frac{c_3 s}{V(x,R) R^{\al} }, $$
where we used VD to obtain the final expression.

Now choose $\delta$ such that $2 \delta^\al = 2^{-\al} \wedge ( C_1/c_3)$.
Let $\kappa\in (0,1]$, and let $s=\kappa (\delta R)^\al$.
If $x,y \in B''=B(x_0, \half \delta R)$, then
$d(x,y) \le \delta R$. So
$$ p_s(x,y) 
\ge \frac{C_1}{V(x,s^{1/\al})} \wedge \frac{C_1 s}{d(x,y)^\al V(x,d(x,y))} 
 \ge \frac{C_1 s }{ (\delta R)^\al V(x,\delta R)} 
\ge \frac{C_1 \kappa }{V(x_0,R)}. $$
So, 
\begin{equation}   p^B_s(x,y) 
 \ge  \frac{C_1 \kappa }{V(x_0,R)} - \frac{c_3 s}{V(x_0,R) R^{\al} } 
 =  \frac{\kappa }{V(x_0, R)}( C_1 - c_3 \delta^\al)
  \ge \frac{C_1\kappa }{2 V(x_0,R)}.  \label{2.3} 
  \end{equation}

Now let $x_1, y_1 \in B'$. Choose $n= 1+ \lfloor 4/\delta \rfloor$, 
and let $x_1=z_0, z_1 \dots, z_n=y_1$ be a sequence of points in $B'$
with $d(z_{i-1},z_i) = \frac14 \delta R$. Let $B_i=B(z_i, \delta R/4)$,
and note that (\ref{2.3}) implies that 
$p^B_s(x,y) \ge c_4 \kappa / V(x_0,R)$ for $x \in B_{i-1}$, $y \in B_i$. 
A standard chaining argument then gives
\begin{equation} 
p^B_{ns}(x_1,y_1) \ge c_5\kappa^n /V(x_0,R). \label{2.4}
\end{equation}
We have
\begin{equation} 
ns =(1+ \lfloor 4/\delta \rfloor)\kappa (\delta R)^\al
\label{2.5}
\end{equation}
so choosing $\kappa$ such that $ns =\lam R^\al$ we obtain (\ref{2.1}). \qed

\begin{rem}{\rm  As already mentioned in the introduction, the need for $\lam$ 
(and so for $\kappa$) 
only arises when $\al \le 1$; when $\al>1$ the usual chaining argument
with sufficiently small balls allows one to bound $p^B_{\lam
R^\al}(x,y)$ from below once one has the bound (\ref{2.3}) with 
$\kappa=1$. }
\end{rem}

\section{Parabolic Harnack inequality}

In this section we 
show (under the assumption $V(d)$) that PI$(\al)$, UJ$(\al)$
and UJS together imply PHI$(\al)$.

\begin{lem}\label{Lem5.1} Let $0<\al<2$.
Suppose  $V(d)$, PI$(\al)$ and  UJ$(\al)$ 
hold. Then the upper bound UHKP$(\al)$ holds.
\end{lem}

\proof It is  well known  (see for example \cite{SC1})
that PI$(\al)$ implies the Nash inequality
\begin{equation} \label{nash}
 ||f||_2^{2+(2\al/d)}  \le C_N \sE(f,f) ||f||_1^{2\al/d}.  
\end{equation}
Given (\ref{nash}), 
we have UHKP$(\al)$ by the arguments in \cite{BL,CK}.
(See also \cite{BGK,CK2} for a simpler version of the proof.) \qed

\ms Given PI$(\al)$ and VD the argument of \cite{SS} gives a weighted
Poincar\'e inequality.
This takes the following form. Let $x_0\in G$, $R\geq 1$, $B=B(x_0,R)$,
and $$ \varphi_R(x)=c_1(R-d(x,x_0))^+,$$
where
$c_1$ is chosen so that $\sum_{x\in B} \varphi_R(x)=1$. Set 
$$ \ol f_{\varphi_R}  = \sum_{x\in B} f(x) \varphi_R(x) \mu_x. $$ 
Then there exists a constant $C$ not depending
on $R,f$, or $x_0$ such that
\begin{equation}
\sum_{x\in B} |f(x)-\ol f_{\varphi_R}|^2 \mu_x \leq
C\sum_{x,y\in B} (f(x)-f(y))^2 
\big(\varphi_R(x)\land \varphi_R(y)\big) J(x,y). 
\label{3.11A}\end{equation}

\begin{lem}\label{Lem5.2} Suppose  $V(d)$, PI$(\al)$ and  UJ$(\al)$ 
hold. Then 
NDLB$(\al)$ holds.
\end{lem}
\proof This follows  from the weighted
Poincar\'e inequality by a standard argument; see, for example,
\cite{BBCK}, Section 3.
So we have 
\begin{equation}  
p^B_T(x',y') \ge \frac{c}{V(x_0,R/2)},  \q   x',y' \in B(x_0,R/4),~~
T\asymp R^\al. 
\label{5.2}\end{equation}
\qed

\begin{pro}\label{phiproof} Suppose VD, UHKP$(\al)$, 
NDLB$(\al)$ and UJS hold. Then PHI$(\al)$ holds.
\end{pro}

\proof Let $\lam \in (0,1]$, $R\ge 1$, $T= \lam R^\al$, $x_0\in G$, and write:
$$ B_0=B(x_0,R/2), \q B'=B(x_0, 3R/4), \q B = B(x_0,R), $$
and 
$$ Q=Q(x_0,T,R)=[0,T]\times B, \q E=(0,T]\times B'. $$
We consider the space time process on $\bR \times G$
given by $Z_t=(V_0-t, X_t)$, for $t\ge 0$. 

Let $u(t,x)$ be non-negative and caloric on $Q$.
Define the r\'eduite $u_E$ by
$$ u_E(t,x)= \bE^x\big( u(t-T_E, X_{T_E}); T_E< \tau_Q \big), $$
where $T_E$ is the hitting time of $E$ by $Z$, and $\tau_Q$
the exit time by $Z$ from $Q$.
Then $u_E=u$ on $E$, $u_E=0$ on $Q^c$, and $u_E\le u$ on $Q-E$. 

The process $Z$ has as a dual the process $\widehat Z_t=(V_0+t, X_t)$; 
we may therefore apply the results of Chapter VI of \cite{BG}.
The balayage formula gives 
$$ u_E(t,x)=  \int_E p^B_{t-r}(x,y) \nu_E(dr,dy), \q (t,x) \in Q, $$
where $\nu_E$ is a measure on $\overline E$.
We write
$$ \nu_E(dr,dy) = \sum_{z\in B'} \nu_E(dr,z) \delta_z(dy) \mu_z . $$
We can divide each of the measures $ \nu_E(dr,z)$ 
into two parts: an atom at 0, and the remainder.
Given this we can write
\begin{equation}  
 u_E(t,x) = \sum_{z \in B'} p^B_t(x,z) u(0,z) \mu_z 
 +  \sum_{z \in B'} \int_{(0,t]} p^B_{t-r}(x,z) \mu_z \nu_E(dr,z). \label{5.3}
 \end{equation}
To identify $\nu_E(dr,z)$ note that if $(t,x) \in E$
then
\begin{equation} \label{nue}
\frac{\pd u_E}{\pd t} =  \frac{\pd u}{\pd t}
 = \sL  u = \sL (u-u_E) + \sL u_E. 
\end{equation}
In particular $u_E(t,x)$ is differentiable. 
Differentiating (\ref{5.3}) we deduce that 
each  measure $\nu_E(dr,z)$ is absolutely continuous with respect to Lebesgue 
measure, and that, writing \break
$\nu_E(dr,z)= v(r,z)\, dr$,
\begin{equation} 
\frac{\pd u_E}{\pd t}(t,x) = v(t,x) + \sL u_E(t,x). \label{5.4}
\end{equation}
Using (\ref{nue}) this gives
\begin{equation} 
v(t,x) =  \sL (u-u_E)(t,x) 
= \mu_x^{-1} \sum_{z\in B-B'}J(x,z) ( u(t,z)-u_E(t,z)). \label{5.5} 
\end{equation}

Now let $(t_1,x_1) \in Q_-$ and $(t_2,x_2) \in Q_+$. 
To prove the parabolic Harnack inequality it is enough, using (\ref{5.3}), 
to show that: 
\begin{align}
  p^B_{t_1}(x_1,z) 
    &\le C p^B_{t_2}(x_2,z) \q \hbox{ for } z \in B', \label{5.6}\\
  \sum_{z \in B'} p^B_{t_1-r}(x_1,z) v(r,z) \mu_z &\le
  C \sum_{z \in B'} p^B_{t_2-r}(x_2,z) v(r,z) \mu_z,
 \q 0\le r \le t_1. \label{5.7} \end{align}

Of these (\ref{5.6}) is immediate from UHKP$(\al)$ and 
NDLB$(\al)$.
So we consider (\ref{5.7}).
Now, writing $w_t(x)=u(t,x)-u_E(t,x) \ge 0$, $s=t-r$,
and $J w_r(z)=\sum_{y\in B-B'} J(z,y)w_r(y)$,
$$  \sum_{z \in B'} p^B_s(x,z) v(r,z)\mu_z
 = \sum_{z \in B'} \sum_{y \in B-B'} p^B_s(x,z) J(z,y)\mu_z w_r(y)
  = \sum_{z \in B'} p^B_s(x,z) Jw_r (z)\mu_z . $$
Then since $t_2-r\ge t_2-t_1\ge T/4$, using 
NDLB$(\al)$,
\begin{equation} 
\sum_{z \in B'} p^B_{t_2-r}(x,z) Jw_r (z) \mu_z
 \ge c  V^{-1} \sum_{z \in B'} Jw_r(z)\mu_z, \label{5.8} 
 \end{equation}
where $V= V(x_0,R)$. 
To complete the proof of (\ref{5.7}) we need,
writing $s=t_1-r \in [0,T/2]$, to show
\begin{equation} 
\sum_{z \in B'} p^B_s(x,z) Jw_r (z) \mu_z
 \le c  V^{-1} \sum_{z \in B'} Jw_r(z) \mu_z. \label{5.9} 
 \end{equation}

If $s \ge T/8$, then using the upper bound on $p^B$
we obtain (\ref{5.9}). So suppose $s\le T/8$.
Let $B_1=B(x_0,5R/8)$. Then 
\begin{equation} 
\sum_{z \in B'} p^B_s(x,z) Jw_r (z)\mu_z =
 \sum_{z \in B_1} p^B_s(x,z) Jw_r (z)\mu_z  +
 \sum_{z \in B'-B_1} p^B_s(x,z) Jw_r (z)\mu_z. 
 \label{5.10}\end{equation}
 If $z \in B'-B_1$ then $d(x,z) \ge R/8$ and so
$$ p^B_s(x,z) \le \frac{cs}{(R/8)^\al V(x,R/8)}
 \le c' V^{-1}. $$
Hence
\begin{equation}  \sum_{z \in B'-B_1} p^B_s(x,z) Jw_r (z) \mu_z
 \le c V^{-1} \sum_{z \in B'-B_1} Jw_r (z)\mu_z
 \le c V^{-1} \sum_{z \in B'} Jw_r (z)\mu_z. \label{5.11}
 \end{equation}
If $z \in B_1$ then using UJS 
\begin{align*}
 Jw_r(z) &= \sum_{y \in B-B'} J(z,y) w_r(y)  \\
 &\le  \sum_{y \in B-B'} 
  V(z,R/16)^{-1} \sum_{z'\in B(z,R/8)} J(z',y) w_r(y)  \\
 &=  V(z,R/16)^{-1} \sum_{z'\in B(z,R/8)} Jw_r(z') 
 \le cV^{-1}  \sum_{z' \in B'} Jw_r (z'). 
\end{align*}
So, using (\ref{mucond}), 
\begin{equation}  \sum_{z \in B_1} p^B_s(x,z) Jw_r (z)
\mu_z \le  cV^{-1}  \sum_{z' \in B'} Jw_r (z')\mu_{z'}. \label{5.12}
 \end{equation}
Combining (\ref{5.11}) and (\ref{5.12}) proves (\ref{5.9}), 
and hence (\ref{5.7}). \qed

\section{Consequences of the parabolic Harnack inequality }
Throughout this section we assume PHI$(\al)$.

\begin{lem}\label{Lem3.1}
Let $G$ satisfy PHI$(\al)$. Then
the on-diagonal upper bound UHD($\al$) holds:
\begin{equation} 
p_t(x,x) \le \frac{C}{V(x,t^{1/\al})}, \q x \in G,\, 
t > 0.
\label{3.1}\end{equation}
\end{lem}
\proof Let $r=t^{1/\al}$ and 
$\lam=1$. 
Let $u(t,y)=p_t(x,y)$ and apply
PHI$(\al)$ to $u$ in 
$Q=(0,4t)\times B(x,2r)$; this gives 
$$ p_{t}(x,x) \le \sup_{Q_-} u \le 
c_0 \inf_{Q_+} u \le c_0p_{3t}(x,y), \q y \in B(x,r). $$
Integrating over $B=B(x,r)$,
$$ \mu(B) p_{t}(x,x) 
 \le c_0 \sum_{y \in B} p_{3t}(x,y)\mu_y \le c_0. $$
\qed

\ms Let $B\subset G$, and $p^B_t(x,y)$ be the heat kernel for $X$ killed 
on exiting $B$. A key consequence of the parabolic Harnack inequality is
a lower bound for $p^B_t(x,y)$. For continuous processes a standard
argument (see \cite{SC2}, p. 153) is to apply the parabolic Harnack 
inequality to the function
\begin{equation} 
v(s,x) =
\begin{cases}  
\psi(x), &\hbox{ if } s < t/2,  \\
 P_{s-t/2} \psi(x), &\hbox{ if } t/2\le s. 
\end{cases}
\end{equation}
where $\psi=1$ on a ball $B$ and $\psi=0$ on $G-B^*$, where $B^*$
is the ball with the same center as $B$ but radius twice as
large. However,
for $v$ to be caloric one needs $\sL \psi=0$ on $B$, and this 
fails if the process can jump from $B$ to $G-B^*$. 
Instead we use the argument below.

\begin{theo}\label{Thm3.2} Let $\Gam$ satisfy PHI($\al$). 
Then if $x_0\in G$, $T=R^\al$, $B=B(x_0,R)$, 
\begin{equation} 
p^B_T(x',y') \ge \frac{c}{V(x_0,R/2)},  \q   x',y' \in B(x_0,R/2). 
\label{3.2}\end{equation}
\end{theo}

\proof Let $R_0=R/2$ and 
$c_0= 1-(3/4)^{1/\al}$. Set  $r_k = c_0 R_0 (3/4)^{k/\al}$, and let
$$ R_k = R_{k-1}-r_{k-1} = R_0 -\sum_0^{k-1} r_i = R_0 (3/4)^{k/\al}. $$
Let $t_k =r_k^\al$, and let $B_k=B(x_0,R_k)$ for $0\le k < \infty$;
for large $k$ the ball $B_k$ is just $\{x_0\}$. 

Set 
$$ u_n(t,x) =\bP^x( \tau_{B_n}>t) =\sum_{y\in B_n} p^{B_n}_t(x,y) \mu_y , $$
$$ \th_n(t) = \sup_{x\in B_n} u_n(t,x). $$
Since $u_n$ is a sum of caloric functions, $u_n$ is caloric
in $(0,\infty)\times B_n$.
Note that $u_{n+1} \le u_n$ and that $u_n(x,t)$ and $\th(t)$ are
decreasing in $t$. Also note that 
\begin{equation} 
u_n(t,x_0) \ge \bP^{x_0}( X_s =x_0, 0\le s \le t) \ge e^{-J(x_0,G)t}.
\label{3.3}\end{equation}

For any ball $B_j$
\begin{align*}
 u_j(t+s,x) &= \sum_{y\in B_j} p^{B_j}_{t+s}(x,y)\mu_y  
=\sum_{z\in B_j} \sum_{y\in B_j}  p^{B_j}_t(x,z) p^{B_j}_s(z,y)\mu_y \mu_z \\
&= \sum_{z\in B_j}  p^{B_j}_t(x,z) u_j(s,z)\mu_z \le u_j(t,x) \th_j(s). 
\end{align*}
Therefore $ \th_j(t+s) \le \th_j(t) \th_j(s)$. 

Now let $n \ge 0$ and let
$x\in B_{n+1}$. Then $B(x,r_n) \subset B_n$, so that 
$u_n(t,x)$ is caloric in $Q=(0,t_n)\times B(x,r_n)$.
Applying the parabolic Harnack inequality we obtain
$$ u_n(t_n/4,x) \le \sup_{Q_-} u_n \le 
  C_1 \inf_{Q_+} u_n \le C_1 u_n(t_n,x) 
\le C_1\th_n(t_n) \le C_1 \th_n(t_n/3)^3.  $$
Since $t_{n+1}/3= t_n/4$, 
\begin{equation} 
\th_{n+1}(t_{n+1}/3) = \sup_{B_{n+1}}u_{n+1}(t_n/4,x) 
\le \sup_{B_{n}}u_{n}(t_n/4,x)
 \le  C_1 \th_n(t_n/3)^3.  \label{3.4}\end{equation}

Write $a_n= \th_n(t_n/3)$; we have 
\begin{equation} 
a_{n+1} \le C_1 a_n^3, \qq n \ge 0. \label{3.5}\end{equation}
Note that $a_n \le 1$ for all $n$.
Suppose that $a_0\le (C_1 \vee e)^{-1}$. Then
$a_1 \le (C_1 a_0)a_0^2 \le a_0$, and so iterating we deduce
that $C_1 a_n\le 1$ for all $n$. Therefore, using (\ref{3.5}) again,
$$  a_n \le (C_1 a_{n-1}) a_{n-1}^2 \le  a_{n-1}^2 
\le  a_{n-2}^4 \le (a_0)^{2^{n}} \le e^{-2^{n}}. $$ 
Hence
$u_n(t_n, x_0) \le \exp(-2^{n})$ for all $n \ge 0$, 
which contradicts (\ref{3.3}).

So $a_0 \ge  (C_1 \vee e)^{-1} =c_2$, and thus 
$\theta_0(s)\ge c_2$ for $s \in [0,t_0/3]$. 
Let $s = t_0/3 \wedge (T/8)$. Then 
there exists $x' \in B_0=B(x_0,R/2)$ such that 
$$ u_0(s,x') \ge c_2. $$ 
Applying the parabolic Harnack inequality to $u_0$ enough times 
to compare 
$u_0(s,x')$ with $u_0(T/4,x')$ it follows that $u_0(T/4,x') \ge c_4$.
Thus as 
$$ u_0(T/4,x') =\sum_{y \in B_0} p^{B_0}_{T/4}(x',y)\mu_y, $$
writing $V_0=\mu(B_0)$, there exists $y'\in B_0$ with 
$$ c_4 V_0^{-1} \le  p^{B_0}_{T/4}(x',y') \le  p^{B}_{T/4}(x',y'). $$
Applying the parabolic Harnack inequality to 
$v(t,y)=p^{B}_{t}(x',y)$  
in the region 
$(0,T) \times B(x_0,R)$
we obtain
$$ c_4 V_0^{-1} \le  p^{B}_{T/4}(x',y) 
 \le \sup_{Q_-} v \le C_1 \inf_{Q_+} v \le C_1 p^B_{3T/4}(x',y), \, 
\hbox{ for all $y\in B_0$.} $$
Fix $y \in B_0$; applying the parabolic Harnack inequality again to 
$w(t,x)= p^B_{t+ T/4}(x,y)$ in the region $(0,T)\times B(x_0,R)$
we obtain, for any $x\in B_0$,
$$ c_4 (C_1 V_0)^{-1} \le  p^B_{3T/4}(x',y) 
 = w(T/2,x') \le \inf_{Q_-} w \le C_1  \inf_{Q_+} w
 \le C_1 w(3T/4,x) =  C_1 p^B_{T}(x,y),  $$
which completes the proof of (\ref{3.2}). \qed

\begin{cor}\label{theo:Corollary3.3}  Suppose $(\Gam,J)$ satisfies PHI$(\al)$.
Then $\Gam$ satisfies VD.
\end{cor}
\proof This is immediate given (\ref{3.1}) and (\ref{3.2}). Let $R>0$,  
$T=R^\al$ and $x\in G$. Then
$$ \frac{c_1}{V(x,R/2)}\le p^B_T(x,x) \le p_T(x,x)\le 
\frac{c_2}{V(x,R)}, $$
giving VD. \qed

\begin{cor}\label{theo:Corollary3.4}  PHI$(\al)$
implies SB($\alpha$). 
\end{cor}

\proof First, note that SB($\alpha$) is equivalent to the following: 
There exist $c_1,c_2,c_3>0$ such that 
for any $x_0\in G$, $r>0$, and writing  $B=B(x_0,r)$,
\begin{equation}\label{ucequi}
\sup_{x,y\in B}p_t^B(x,y) \le c_1 V(x_0,t^{1/\alpha})^{-1}, 
\q \hbox{ for all } c_2r^\al \le t \le c_3 r^\al.
\end{equation} 
Now $p_t(x,y)^2 \le p_t(x,x)p_t(y,y)$. So, by Lemma \ref{Lem3.1},
for $x, y \in B$,
\begin{equation*}
  p_t^B(x,y)^2 \le 
  p_t(x,y)^2 \le \frac{c }{V(x,t^{1/\al})V(y,t^{1/\al})  }
 \le  \frac{c_4 }{V(x_0,t^{1/\al})^2 },
\end{equation*}
where we used (\ref{vd2}) in the last line. \qed

\begin{lem}\label{Lem3.4} Suppose $G$  satisfies PHI$(\al)$.
Then
\begin{equation} 
c_1 R^\al \le \bE^x \tau_{B(x,R)} \le c_2 R^{\al}, \qq R\geq 1. \label{3.6}
\end{equation}
\end{lem}
\proof Let $B=B(x_0,R)$ and $B'=B(x_0,R/2)$; then if $T=(2R)^\al$, 
(\ref{3.2}) gives
$$ p_T(x,y) \ge c V(x_0,R)^{-1}, \q x,y \in B. $$
Fix $y_0$ with $d(x_0,y_0)=[3R/4]$; then if $x\in B$
$$ \bP^x( X_T \not\in B') \ge \bP^x(X_T \in B(y_0,R/4)) 
= \sum_{y\in B(y_0,R/4)} p_T(x,y)\mu_y 
\ge c\frac{ V(y_0,R)}{V(x_0,R)}\ge c_3. $$
So  we have
$\bP^x( \tau_{B'} > T) \le 1-c_3$ for all $x \in B'$. Hence
by the Markov property
$\bP^x( \tau_{B'} > kT) \le (1-c_3)^k$, and thus $\bE^x \tau_{B'} \le c_4 T$.
Since $\bE^{x_0} \tau_{B(x_0,R/2)} = \bE^{x_0} \tau_{B'}$, replacing
$R/2$ by $R$ gives the upper bound in (\ref{3.6}).

The lower bound is easy; Theorem \ref{Thm3.2} gives 
$$   \bP^x( \tau_{B(x,R)} > R^\al) 
 = \sum_{y\in B(x,R)} p^B_{R^\al}(x,y)\mu_y
\ge c_4 >0, $$
and thus $\bE^x \tau_{B(x,R)} \ge c_4 R^\al$. \qed

\begin{rem}{\rm Lemma \ref{Lem3.4} shows that $(\Gam,J)$ 
cannot satisfy PHI$(\al)$ for two different values of $\al$.}
\end{rem}

\begin{pro}\label{theo:Proposition3.6} 
Suppose $(\Gam,J)$ satisfies PHI$(\al)$. Then UJS holds.
\end{pro}

\proof 
Let $A \subset G$ and $f(t,x)$, $t\in \bR_+$, $x \in G-A$,
be a  bounded non-negative function. Consider the equations
\begin{align}
 \frac{\pd u}{\pd t}(t,x) &= \sL u(t,x), \qq x \in A,  \label{3.8}\\
 u(0,x)&= 0,\qq x \in A, \label{3.9}\\
 u(t,x) &= f(t,x), \qq x \in A^c. \label{3.10}
 \end{align}
Then $u$ is caloric on $(0,\infty) \times A$ and 
\begin{equation} 
u(t,x) =\bE^x(f(t-\tau_A,X_{\tau_A});\tau_A \le t). \label{3.11}
\end{equation}

Let $x_0, y_0 \in G$ and $R\le d(x_0,y_0)$. Take
$A=B(x_0,R)$, let $T=R^\al$, $h>0$ be small and define $f_h(t,x)$ by
$$  f_h(t,x) = 1_{(x=y_0)} 1_{(T/2-h,T/2)}(t), \q x \in G-B. $$
Let $u_h(t,x)$ be the solution of (\ref{3.8})--(\ref{3.10}). Thus
$$ u_h(t,x)= \bP^x( X_{\tau_B} =y_0, t-T/2< \tau_B < t+h-T/2). $$
Since $C_J < \infty$ we have
\begin{equation} 
  \lim_{ h \downarrow 0} h^{-1} u_h(T/2,x) = \mu_x^{-1} J(x,y_0). \label{3.12}
\end{equation}
Applying the parabolic Harnack inequality to $u$ in 
$(0,T)\times B(x_0, R)$ we obtain
$$  u_h(T/2,x_0) \le C_1 u_h(T,x_0). $$
Now by (\ref{ucequi})
$$ u_h(T,x_0) = \sum_{z\in B} p^B_{T/2}(x_0,z) u_h(T/2,z) \mu_z
  \le c \mu(B)^{-1}  \sum_{z\in B} u_h(T/2,z)\mu_z . $$
Thus
$$  u_h(T/2,x_0) \le c \mu(B)^{-1}  \sum_{z\in B} u_h(T/2,z)\mu_z, $$
and using (\ref{3.12}) gives
\begin{equation}
 J(x_0,y_0) \le \frac{c\mu_{x_0}}{V(x_0,R)}\sum_{z \in B(x_0,R)} J(z,y), 
\end{equation}
proving UJS. \qed

\begin{lem}\label{Lem3.7} Suppose PHI($\al$) holds.
Let $B=B(x_0,R)$, and $B'=B(x_0,R/2)$.  Then
\begin{equation} 
\sum_{y \in B'} J(y,G-B) \le  \frac{c \mu(B')}{R^{\al}} . \label{3.14}
\end{equation}
\end{lem}

\proof Let $\tau=\tau_B$, and consider the martingale
$$ M_t = 1_{[\tau,\infty)}(t) -\int_0^t 1_{(\tau> s)} \mu_{X_s}^{-1} J(X_s,G-B) ds. $$
Then $\bE^x M_t=0$, and hence
\begin{equation} 
1 \ge \bE^x \int_0^t 1_{(\tau> s)} \mu_{X_s}^{-1} J(X_s,G-B)ds
 = \int_0^t \sum_{y \in B} p^B_s(x,y) J(y,G-B)ds  . \label{3.15}
 \end{equation}

Using the lower bound (\ref{3.2}), and writing $T=R^\al$, 
$$ 1\ge \sum_{y \in B'} J(y,G-B) \int_{T/2}^T   p^B_s(x_0,y) ds
 \ge c T \mu(B')^{-1} \sum_{y\in B'} J(y,G-B) . $$
\qed

\begin{pro}\label{theo:Proposition3.8} If PHI$(\al)$ holds then 
UJ$(\al)$ holds, i.e.,
\begin{equation} \label{jub4}
 J(x,y) \le \frac{c\mu_x \mu_y } {d(x,y)^\al V(x,d(x,y))}. 
\end{equation}
\end{pro}

\proof  
Using (\ref{1.1}) and (\ref{mucond}), (\ref{jub4}) holds if 
$d(x,y) \le 3$. If $d(x,y) >3$ let  $r=\lfloor d(x,y)/3 \rfloor$. 
Then using Proposition \ref{theo:Proposition3.6} twice
and (\ref{3.14}) once,
\begin{align*}
J(x,y) &\le \frac{c \mu_{x}}{V(x,r)} \sum_{x' \in B(x,r)} J(x',y) \\
 &\le \frac{c \mu_{x}}{V(x,r)} \sum_{x' \in B(x,r)} 
 \frac{c \mu_{y}}{V(y,r)} \sum_{y' \in B(y,r)} J(x',y') \\
 &\le \frac{c \mu_x \mu_y }{V(x,r)V(y,r)}  
  \sum_{x' \in B(x,r)} J(x',B(x_0, 2r)^c)  \\
 &\le  \frac{c \mu_x \mu_y }{V(x,r)V(y,r)} \frac{V(x,r)}{r^\al} 
 =  \frac{c \mu_x \mu_y }{r^\al V(y,r)}.
\end{align*}
Using (\ref{vd2}) completes the proof. \qed

\begin{lem}\label{Lem3.9}
Suppose PHI$(\al)$ holds. Then PI$(\al)$ holds.
\end{lem}

\proof Let $B$ be a ball and let $\ol X$ denote the process $X$
`reflected on the boundary of $B$.' That is, $\ol X$ is the 
process with jump rates 
\begin{equation} 
 \ol J(x,y)= 
\begin{cases}  
 J(x,y), &\hbox{ if } x,y \in B, \\ 
   0,     &\hbox{ otherwise}. \\
\end{cases}
\end{equation}
Write $\ol p_t(x,y)$ for the heat kernel of $\ol X$. Then
by Theorem \ref{Thm3.2}  
$$ \ol p_t(x,y) \ge p^B_t(x,y) \ge 
 \frac{c}{V(x_0,R/2)},  \q   x,y \in B(x_0,R/2). 
\eqno(3.16) $$
This lower bound then  gives
PI$(\al)$ by a standard argument, as in  \cite{SC2}, p. 159--160.
\qed

\sms
We summarize the results of this section in the following
Proposition.

\begin{pro} \label{phicons}
Suppose PHI($\al$) holds. 
\newline (a) $(\Gam,J)$ satisfies VD, UHD$(\al)$,
NDLB$(\al)$, SB$(\al)$, $E_\al$, UJS, UJ$(\al)$ and PI$(\al)$.
\newline (b) Suppose that in addition  $(\Gam,J)$ satisfies
$V(d)$. Then UHKP$(\al)$ holds.
\end{pro}

\section{Consequences of the on-diagonal upper bound}

In this section, we assume that $(\Gam,J)$ satisfies 
VD and $SB(\alpha)$. 
For $A \subset G$ let
\begin{equation} 
 \sF_A =\{u\in L^2(G,\mu): u|_{G-A}=0\}. 
 \end{equation}

\begin{lem}\label{Lem4.0} 
Suppose VD and $SB(\al)$ hold. 
Then there exists $c_1>0$ such that for all $x_0 \in G$, $r \ge 1$,
\begin{equation}\label{LNal}
\sE (u,u)
\ge c_1 \frac {\|u\|_2^2}{r^\al},\q  \mbox{for all }~u\in \sF_{B(x_0,r)}.
\end{equation}
\end{lem}

\proof Let $\lam\ge 1$ (to be chosen later) and $B=B(x_0,\lam r)$. 
Using the log-convexity of $t\mapsto \|P^B_tu\|_2^2$
(see \cite{Cou} Lemma 3.2 for the proof), we have 
\begin{equation} 
\frac{(P^B_tu,u)}{\|u\|_2^2}
 \ge \exp \Big(-\frac{\sE (u,u)}{\|u\|_2^2}t\Big),\qquad 
 \hbox{ for all } u\in \sF_B, \q t \ge 0.\label{4.1}
 \end{equation}
By interpolating $SB(\alpha)$ with $\|P^B_t \|_{1\to 1}\le 1$, and 
using VD, we obtain 
$$\|P^B_tu\|_2^2 \le 
c_1 V(x_0,\lam r)^{-1}\|u\|_1^2,
 \qquad c_2(\lam r)^\al\le t \le c_3 (\lam r)^\al.$$
Substituting this into (\ref{4.1}) with $2t$ instead of $t$,  
\begin{equation}  \frac{c_1 \|u\|_1^2} {V(x_0,\lam r)\|u\|_2^2}
 \ge \exp \Big(-2\frac{\sE (u,u)}{\|u\|_2^2}t\Big),
\qquad c_2 (\lam r)^\al \le t \le c_3 (\lam r)^\al,
~u\in \sF_B\cap L^1.\label{subst}\end{equation}
Let $u\in \sF_{B(x_0,r)}\cap L^1$, and $t=c_3 (\lam r)^\al$; then
using the Cauchy-Schwarz inequality, 
\begin{equation}\label{vratio}
\frac{c_1 \|u\|_1^2}{V(x_0, \lam r) \|u\|_2^2}  \le 
\frac{c_1 V(x_0,r)}{V(x_0,\lam r)}. 
\end{equation}
As $r\ge 1$, using (\ref{vd3}),
we can choose $\lam$ so that the right hand side of
(\ref{vratio}) is less than $e^{-1}$. 
(\ref{subst}) with $t= c_3 (\lam r)^\al$ and (\ref{vratio}) 
now give (\ref{LNal}). \qed

\ms
Let 
\begin{equation*}
  M_1(x,r) = \sum_{y\in B(x,r)} d(x,y)^2 J(x,y),  \qquad
  M_2(x,r) =\sum_{y\in B(x,r)^c} J(x,y).  
  \end{equation*}

\begin{lem}\label{Lem4.1} Suppose VD and $SB(\al)$ hold. 
Then for all $x_0$, $r$
\begin{equation} 
\frac{c_1 V(x_0,r)}{r^{\al}} \le  \sum_{x\in B(x_0,r)} 
 (r^{-2}M_1(x,r)+M_2(x,r)) . \label{4.3}
 \end{equation}
\end{lem}

\proof 
Let $x_0\in G$, and let $r>0$. Consider the function 
$$ f(y) = (1 - r^{-1} d(x_0,y))_+. $$
Let $A(x)=\{y : d(x_0,y) \ge d(x_0,x)\}$, and
$$  \Gam f(x)=  \sum_{y\in A(x)} (f(x)-f(y))^2 J(x,y). $$
Then $f(x)=0$ and $\Gam f(x)=0$ if $x \in B(x_0,r)^c$, so 
$$ \sE (f,f) \le 2 \sum_{x \in G}  \sum_{y\in A(x)} (f(x)-f(y))^2 J(x,y)
 = \sum_{x \in B(x_0,r)} \Gam f(x) . $$
Since $|f(x)-f(y)| \le (c_0r)^{-1} d(x,y)$, and $0\le f \le 1$,
for $x\in B(x_0,r)$ we have 
$$  \Gam f(x) \le c\sum_{B(x,r)} (d(x,y)/r)^2 J(x,y)
 +  \sum_{B(x,r)^c}  J(x,y)\le c r^{-2}M_1(x,r) + M_2(x,r). $$
Combining these inequalities 
\begin{equation} \sE (f,f) 
 \le C \sum_{B(x_0,r)} (r^{-2}M_1(x,r)+M_2(x,r)) , \label{4.2}
 \end{equation}
and (\ref{4.3}) follows by Lemma \ref{Lem4.0}.

\begin{pro}\label{theo:Proposition4.3} 
Let $0<\alpha<2$ and assume VD. Suppose $SB(\alpha)$ and UJ$(\al)$ hold. 
Then there exist $\delta>0$,  $\lam<\infty$ (depending only on $\al$ and on the
constants $C$ in $SB(\al)$ and UJ$(\al)$) so that for 
all $x_0 \in G$, $r \ge 1$,
\begin{equation} 
\sum_{x\in B(x_0,r)} \Big[ \sum_{y\in B(x,\lam r)-B(x,\delta r)} J(x,y)  \Big]
\ge c_1  \frac {V(x_0,r)}{r^\alpha}.   
\label{4.4}
\end{equation}
\end{pro}

\proof Note that the term in brackets on the left side of (\ref{4.4}) is
$M_2(x,\delta r)-M_2(x,\lam r)$.

Using UJ$(\alpha)$ and VD, and the fact that $\al<2$, we have 
\begin{align*}
M_1(x,r) &\le  
c_1 \sum_{ B(x, r)} \frac{d(x,y)^{2-\al}}{ V (x,d(x, y))}\\
&\le  c_2 \sum_{i=0}^\infty 
 \sum_{ B(x,2^{-i} r)-B(x,2^{-i-1} r)}  \frac{d(x,y)^{2-\al}}{ V (x,d(x, y))}\\
&\le c_2  \sum_{i=0}^\infty (2^{-i}r)^{2-\al} 
  \frac{V(x,r 2^{-i})}{ V (x, 2^{-i-1}r )} \\
&\le c_3  \sum_{i=0}^\infty (2^{-i}r)^{2-\al}  = c_4 r^{2-\al},
\end{align*}
and
\begin{align*}
M_2(x,r)
&\le  c_5 \sum_{i=0}^\infty 
 \sum_{ B(x,2^{i+1} r)-B(x,2^{i} r)}  \frac{1} { d(x,y)^{\al} V(x,d(x, y))}\\
&\le  c_5 \sum_{i=0}^\infty (2^i r)^{-\al}
 \frac{ V(x, 2^{i+1} r)}{V(x, 2^i r)} \le c_6 r^{-\alpha}.
\end{align*}
So, for $x\in B(x_0,r)$,
\begin{align*}
 r^{-2} M_1(x,r) + M_2(x,r) 
  &= \sum_y (r^{-2} d(x,y)^2 \wedge 1) J(x,y) \\
&\le r^{-2} M_1(x, \delta r)  + M_2(x,\delta r) \\
&\le r^{-\al} ( c_4 \delta^{2-\al} + c_6 \lam^{-\al})
 + M_2(x,\delta r) -M_2(x,\lam r).
    \end{align*}
Now choose $\delta>0$  small enough and $\lambda>0$ large
enough so that
$$  c_4 \delta^{2-\al} + c_6 \lam^{-\al} 
\le \half C_M^{-1} c_{5.2.1}; $$
then summing over $x\in B(x_0,r)$ and using (\ref{4.3}) we deduce
(\ref{4.4}). \qed

\begin{lem} \label{Jloc}
Suppose VD and LJS hold. 
Then if $x\sim y$ and $z \neq x,y$,
$$  J(x,z) \ge c J(y,z). $$
\end{lem}

\proof If $d(x,z)=1$ then by (\ref{ljs2}) $ J(x,z)\ge c_0$,
while by (\ref{1.1}) $J(y,z)\le C_J$.
If $d(x,z)\ge 2$ then 
by (\ref{ljs1}) and VD
$$  J(x,z) \ge c \frac{\mu_x}{V(x,1)} \sum_{w \in B(x,1)} J(w,z)
 \ge  c \frac{\mu_x}{V(x,1)} J(y,z)\ge c' J(y,z). $$
\qed

\begin{lem}\label{LemJS} Suppose VD, LJS and UJS hold. Then
JS holds.
\end{lem}

\proof We prove (\ref{js}). 
Let $d(x_0,y)=R$. If $R\le 4$ then we can use Lemma \ref{Jloc},
so suppose $R\ge 4$. Suppose first that $d(x_1,x_0)\le R/4$. Then
writing $s=R/4$, and using UJS and LJS, 
$$ J(x_1,y) \le \frac{c \mu_{x_1}}{V(x_1,s)} \sum_{z \in B(x_1,s)} J(z,y)
 \le \frac{c \mu_{x_1}}{V(x_1,s)} \sum_{z \in B(x_0,2s)} J(z,y)
 \le \frac{c \mu_{x_1}V(x_1,2s)}{\mu_{x_0} V(x_1,s)}J(x_0,y). $$
Using VD and (\ref{mucond}) then gives $J(x_1,y) \le c_1 J(x_0,y)$,
proving  (\ref{js}). If  $d(x_1,x_0)\ge R/4$ then
 (\ref{js}) follows by an easy chaining argument. \qed

\ms We need a general lemma on symmetric functions on $G\times G$ 
which satisfy conditions similar to JS. See \cite{Ba} for a
similar argument.

\begin{lem}\label{Lem4.4} Let $g: G \times G \to \bR_+$
satisfy $g(x,y)=g(y,x)$ for all $x,y$, and also the
conditions 
\begin{align} \label{gloc1}
g(x,y) &\ge c_0,  \q \hbox { if } x\sim y, \\
\label{gloc2}
g(x,z) &\ge c_0 g(y,z),  \q \hbox { if } x\sim y, z \neq x, y. 
\end{align}
 Suppose that 
for some $\kappa\in (0,1)$, $c_1 < \infty$, 
\begin{equation} 
g(x,y) \le c_1 g(x,y') \q \hbox { if } \, d(y,y') 
\le  \kappa d(x,y). \label{4.8}\end{equation}
Then given $0<\delta<\lam<\infty$, there exists a constant 
$C_1$, depending only on $c_1, \kappa, \delta, \lam$, such that
the following holds. 
If $x_0, y_0 \in G$ with $d(x_0,y_0)=r$ then
\begin{equation}
C_1^{-1} g(x_0,y_0) \le g(x,y) \le  C_1 g(x_0,y_0) \q 
 \hbox{ whenever } x,y \in B(x_0, \lam r), \, d(x,y) \ge \delta r. 
 \label{4.9}
\end{equation}
\end{lem}

\proof Let $H$ be the metric space obtained by replacing each
edge of the graph $G$ by a line segment of length 1. (In
\cite{BB} this is called the `cable system' of $G$.) We write
$d$ for the metric on $H$. Extend
$g$ to a function $h$ on $H \times H$ by linearity on each cable;
then the conditions (\ref{gloc1}) and  (\ref{gloc2})
imply that (\ref{4.8}) 
also holds for $h$. So it is now enough
to prove the Lemma for $h$.

We can assume $\delta \le \frac14$ and $\lam\ge 2$.
Also, by an easy chaining argument we can assume $\kappa=\frac12$.
Note first that 
(\ref{4.8}) implies
\begin{equation}
 h(x,y) \asymp h(x,y') \q 
 \hbox{ whenever } d(y,y') \le 
 \half (d(x,y)\vee d(x,y')).
 \label{44a}
\end{equation}

Given $x,y \in H$ let $\gamma(x,y)$ denote a shortest (geodesic) 
path between $x$ and $y$. 
Suppose $x, y \in H$, $d(x,y)=s$, and $z \in \gamma(x,y)$ with
$d(x,z)\ge s/2$. Then by (\ref{44a})
\begin{equation}
 h(x,y) \asymp h(x,z). 
\end{equation}
Using this repeatedly, we can compare $h$ on any geodesic path.
More precisely, if $x, y \in H$, $d(x,y)=s$ then we have
\begin{equation}
 h(x,y) \asymp h(x',y') \q 
\hbox{ for $x', y' \in \gamma(x,y)$ with $d(x',y')\ge \half \delta s$.}
\label{4.10}\end{equation}

Now let $x_0, y_0 \in H$ with  $d(x_0,y_0)=r$,
and  $x_1,y_1 \in B(x_0, \lam r)$, $d(y_1,x_1)\ge \delta r$.
As $G$ is infinite there exists $w \in H$ with
$d(x_0,w) = 5 \lam r$.
Suppose we can prove:
\begin{equation} \label{44c}
 h(x',y') \asymp h(x',w) \q 
\hbox{ for all  $x',y' \in B(x_0,\lam r)$ with } d(x',y') \ge \delta r. 
\end{equation}
Then we have $h(x_j,y_j) \asymp h(x_j,w)$ for $j=0,1$. But since
$d(x_0,x_1) \le \lam r \le d(x_0,w)$, using (\ref{44a}) we have
$h(x_0,w) \asymp h(x_1,w)$, and so we obtain
\begin{equation}
 h(x_0,y_0) \asymp h(x_1,y_1). 
\end{equation}

It remains to prove (\ref{44c}).
Suppose first that 
\begin{equation} \label{44ass}
  d(x',z) \ge \delta r/2 \q \hbox { for all } z \in \gam(y',w).
\end{equation}
Then chaining the relation (\ref{44a}) along $\gam(y',w)$
gives $h(x',y') \asymp h(x',w)$, proving (\ref{44c}).

Now suppose that (\ref{44ass}) fails. Then there exists 
$z \in \gam(y',w)$ with $d(x',z)\le \delta r/2$. 
By (\ref{44a}) $h(y',x') \asymp h(y',z)$. Also, 
$d(y',z)\ge \delta r/2$, 
so using (\ref{4.10}) we obtain
$h(y',z) \asymp h(y',w)$. Finally, as 
$d(x',y') \le 2 \lam r$ and $d(y',w) \ge 4 \lam r$,
(\ref{44a}) gives $h(x',w)\asymp h(y',w)$.
Combining these inequalities gives (\ref{44c}) in this
case also. \qed

\begin{lem}\label{Lem4.6} 
(a) VD, (\ref{4.4}) and JS imply LJ$(\al)$.
\nl (b) Assume $0<\alpha<2$. Then VD, $SB(\al)$, UJ$(\al)$ and JS imply J$(\al)$.
\end{lem}
\proof (a) Let  $x_0, y_0 \in H$, and  $d(x_0,y_0)=r$. Then
by (\ref{4.4}), JS and Lemma \ref{Lem4.4},
\begin{align*}
  c \frac{V(x_0,r)}{r^\al} &\le
  \sum_{x\in B(x_0,r)} 
   \Big( \sum_{y\in B(x,\lam r)-B(x,\delta r)} J(x,y) \Big)\\
 &\le c' J(x_0,y_0)   \sum_{x\in B(x_0,r)} V(x, \lam r)\\
 &\le  c' J(x_0,y_0) V(x_0,r) V(x_0, (1+\lam)r), 
 \end{align*}
and using VD we obtain LJ$(\al)$. 
\nl (b) By Proposition \ref{theo:Proposition4.3}, VD, $SB(\al)$ and UJ$(\al)$
imply (\ref{4.4}). (a) now gives LJ$(\al)$,  and so J$(\al)$ holds. \qed

\begin{pro}\label{theo:Proposition4.7} Assume $0<\alpha<2$. 
\nl (a) PHI$(\al)$ implies (\ref{4.4}).
\nl (b) PHI$(\al)$ and LJS imply J$(\al)$.
\end{pro}

\proof (a) Assume $ PHI(\al)$. Then by Corollary \ref{theo:Corollary3.3}, 
Corollary \ref{theo:Corollary3.4}
and Proposition \ref{theo:Proposition3.8}, VD, $SB(\al)$ and UJ$(\al)$ hold.
Hence, by Proposition \ref{theo:Proposition4.3}, (\ref{4.4}) holds.
\nl (b) Since  PHI$(\al)$ also implies UJS (due to Proposition
\ref{theo:Proposition3.6}) and VD (due to Corollary
\ref{theo:Corollary3.3}), by Lemma \ref{LemJS}
we obtain JS and hence, by Lemma \ref{Lem4.6}(a), J$(\al)$ holds. \qed

\med {\it Proof of Theorem \ref{Thm1.2}.}
(a) $\Rightarrow$ (b). This has been proved in the
context of Markov processes on $\bZ^d$ and on $d$-sets in \cite{BL,CK,CK2}.
The transfer of these arguments to a graph satisfying $V(d)$ is 
straightforward. 
(b) $\Rightarrow$ (a) is immediate, since we
have $J(x,y) =\mu_x \mu_y \lim_{t \to 0} t^{-1} p_t(x,y)$.

Now suppose  J$(\al)$ and HKP$(\al)$ hold. 
Then UJS holds, and by Lemma \ref{Lem2.1},  
NDLB$(\al)$ holds.
Therefore, by Proposition \ref{phiproof} PHI$(\al)$ holds.
Thus  ((a) and (b)) together imply (c).
Finally, by Proposition \ref{theo:Proposition4.7} 
we have (c)  $\Rightarrow$ (a). 
\qed 

\med {\it Proof of Theorem \ref{Thm1.3}.} 
That (a) implies (b) is immediate from Proposition \ref{phicons}.
We remark that this does not need $V(d)$ or $\al<2$. 

\noindent (b) $\Rightarrow$ (a).
This follows by combining Lemmas \ref{Lem5.1} and 
\ref{Lem5.2}, and 
Proposition \ref{phiproof}. \qed

\bigb {\it Proof of Theorem \ref{Thm1.1}.}
(a) is contained in the implication (b) $\Rightarrow$ (c)
in Theorem  \ref{Thm1.2}. (Note that this part of the argument
does not use $V(d)$ or $\al<2$.)
(b) is immediate from Proposition \ref{phicons}. \qed

\begin{rem} 
{\rm One might ask if the three conditions
in Theorem \ref{Thm1.3} are independent.
\begin{itemize}
\setlength{\itemsep}{-\parsep}
\item[1.] If we drop  UJ$(\al)$ then we have no upper bound on $J$.
If UJS and PI$(\al)$ hold, then since PI$(\al)$ implies PI$(\al')$
for any $\al'>\al$, we have UJS and PI$(\al')$ for all 
$\al'\ge \al$. 
However, by the remark following 
Lemma \ref{Lem3.4} we cannot have PHI$(\al')$ for any $\al'>\al$.
\item[2.] If we drop PI$(\al)$ then we have no lower bound on $J$.
We can set  $J(x,y)=d(x,y)^{-d-\al}$, and note that
UJ$(\al')$ and UJS hold for any $\al'<\al$.
\item[3.]  We do not have an example to prove that UJS
is independent of PI$(\al)$ and UJ$(\al)$.  Note that since PI$(\al)$
implies a Nash inequality, (\ref{5.6}) does hold, and 
so gives some kind of lower
bound on $J$. We `only just' needed to use UJS in the proof 
of Proposition \ref{phiproof}, 
to control $Jw_r(z)$ when $z$ and $y$ are far apart.
\end{itemize}
}
\end{rem}

\begin{rem} 
{\rm In the definition of PHI$(\al)$ we included
a parameter $\lam \in (0,1]$. Suppose we call
PHI$(1,\al)$ the PHI just with $\lam=1$.
Then PHI$(1,\al)$ is enough
to obtain UJ$(\al)$, PI$(\al)$ and UJS, 
so Theorem \ref{Thm1.3} gives that
PHI$(1,\al)$ and  PHI$(\al)$ are equivalent.}
\end{rem}

\section{Counterexamples.}
\med {\sl PHI$(\al)$ does not imply HKP$(\al)$.}

Let $G=\bZ^d$, $\al \in (0,2)$ and let
$J_1(x,y)=|x-y|^{-d-\al}$ for $x\neq y$.
Note that $V(d)$ and J$(\al)$ hold for $J_1$. 
So by Theorem 1.2 we have that HKP$(\al)$ and PHI$(\al)$
hold for $J_1$. (Of course, for this example this was already 
well known.) So, by Theorem \ref{Thm1.3}, 
$PI(\al)$, UJS, and UJ$(\al)$ hold for $J_1$.

Choose $R \in 2\bN$, with $R \gg 1$, 
and let $y_0=(R, 0, \dots ,0)$ be on the $x_1$--axis
with $|y_0|=d(0,y_0)=R$.
Then let 
\begin{equation} 
 J(x,y) = 
\begin{cases}  
 J_1(x,y), &\hbox{ if } \{x,y\} \neq \{0,y_0\}, \\
 0, &\hbox{ if }  \{x,y\}= \{0,y_0\}.
\end{cases}
\end{equation}
(So we just suppress jumps between $0$ and $y_0$.)
Since J$(\al)$ fails for $J$, by Theorem \ref{Thm1.2} 
HKP$(\al)$ must fail.

However, PHI$(\al)$ does hold.
To see this we use Theorem \ref{Thm1.3},
and verify that UJ$(\al)$, UJS and $PI(\al)$ all hold.
First, as $J \le J_1$, UJ$(\al)$ is immediate. 
Since $J(0,\cdot)$ has only been modified from $J_1(0,\cdot)$
at one point, it is straightforward to verify that UJS still holds for $J$.

Finally, to verify $PI(\al)$, let $B=B(x_0,r)$ be a ball in $\bZ^d$,
and $f: B \to \bR$. If $B$ does not contain 
both $0$ and $y_0$ then
$$ \sum_{x\in B} \sum_{y \in B} (f(x)-f(y))^2 J_1(x,y)
 = \sum_{x\in B} \sum_{y \in B} (f(x)-f(y))^2 J(x,y), $$
so the Poincar\'e inequality for $J$ follows from that for $J_1$. Now suppose
that $0, y_0 \in B$. Then let $y_1$ be the mid-point of the
line between $0$ and $y_0$. We have
\begin{align*}
 (f(0)-f(y_0))^2 J_1(0,y_0) 
 &\le 2 ((f(0)-f(y_1))^2+ (f(y_1)-f(y_0))^2) R^{-d-\al}  \\
 &= 2^{d+\al+1}\Big( ((f(0)-f(y_1))^2 J(0,y_1) 
+ (f(y_1)-f(y_0))^2) J(y_1,y_0)\Big)\\
 &\le  2^{d+\al+1} \sum_{x\in B} \sum_{y \in B} (f(x)-f(y))^2 J(x,y).
 \end{align*}
Thus 
$$ \sum_{x\in B} \sum_{y \in B} (f(x)-f(y))^2 J_1(x,y)
 \le c  \sum_{x\in B} \sum_{y \in B} (f(x)-f(y))^2 J(x,y), $$
and this implies that the Poincar\'e inequality holds for $J$.

\med {\sl EHI + $E_\al$ + $V(d)$ does not imply  PHI$(\al)$}.

\smallskip
We only give an outline of this example. 
Let $G= \bZ$, $\al \in (1,2)$ and $J_0(x,y) = |x-y|^{-1-\al}$.
Let $R_1 \gg 1$, and set 
$J_1(x,y)= (\log R_1) R_1^{-1-\al} 1_{(|x-y|=R_1)}$.
Let $X^{(i)}_t$, $i=0,1$ be independent 
processes associated with
the jump kernels $J_i$. Let $X=X^{(0)}+X^{(1)}$; this is the
process with jump kernel $J=J_0+J_1$.
We take $\mu_x=1$ for all $x$. 

We begin by remarking that HKP$(\al)$ does hold for
$X^{(0)}$. In addition this process is `strongly recurrent':
one has
\begin{equation} \label{srec}
  \bP^x( T^{(0)}_y \le \tau^{(0)}_R) \ge p_0>0 
\q \hbox{ for } x,y \in [-R/2,R/2], 
\end{equation}
where $ \tau^{(0)}_R$ is the exit time from 
$B(0,R)=[-R,R]$ by  $X^{(0)}$, and $ T^{(0)}_y$
is the hitting time of $y$ by $X^{(0)}$.

We now show that $X$ satisfies $E_\al$. The upper bound 
is easy. Since $J \ge J_0$, the Nash inequality (\ref{nash})
holds for $X$. Hence, by \cite{CKS}, the transition
density of $X$ satisfies
$$   p_t(x,y) \le c_1 t^{-1/\al},\q  t > 0.  $$
So taking $c_2$ large enough, if $t=c_2 r^\al$ then
$\bP^x( X_t \in B(0,r)) \le \half$ for any $x \in \bZ$,
and the upper bound $E^x \tau_{B(0,r)} \le c_3 r^\al$
follows.

For the lower bound, note that the condition $E_\al$ for
$X^{(0)}$ implies that
$$ \bP^0( \tau^{(0)}_R \le \lam R^\al )\ge c_4 \lam. $$
Thus there exists $c_5>0$ such that, 
$\bP^0(   \tau^{(0)}_R \ge c_5 R^\al )\ge c_5$. 

Let $\delta =  R_1^{-1-\al} \log R_1$. 
By Doob's inequality, writing $Y_t =  \sup_{s \le t}  |X^{(1)}_s|$,
\begin{equation} \label{x1bound}
\bE^0 Y_T^2 \le  4 R_1^2 \delta T, 
\end{equation}
and so 
\begin{equation} \label{ycheb}
\bP^0(  Y_T \ge \lam)
 \le \frac{ 4T \log R_1}{ \lam^2 R_1^{\al-1}}. 
\end{equation}
So 
$$   \bP^0(   \tau_{R}  \ge c_5 R^\al )
 \ge \bP^0(  \tau^{(0)}_{R/2} \ge  c_5 R^\al,  Y_T \le R/2)
 \ge \half c_5, $$
provided $R_1$ is large enough.
This establishes the lower bound in $E_\al$ for $X$.

To prove EHI it is enough to prove (\ref{srec}) for $X$,
and using translation invariance and chaining it is enough
to prove that
\begin{equation} \label{sr2}
  \bP^x(  T_0 \le  \tau_R) \ge p_1>0
\q \hbox{ for }  x \in [-R/4, R/4]. 
\end{equation}

As the whole argument is more lengthy than this 
counterexample deserves, we only sketch the main ideas.
We note that there exists $\theta \in (0,1)$ such that
\begin{equation} \label{potx0}
  \bP^x(  T^{(0)}_0 > \tau^{(0)}_R) \le c_6( |x|/r)^\theta,
\q \hbox{ for }  x \in [-R/4, R/4]. 
\end{equation}

Fix an interval $R$, and let $x \in [-R/2,R/2]$. 
We concentrate on the case when $R_1 \ll |x| \ll R$.
Choose $\eps>0$ small, and let 
$r=r(x)=|x|^{1+\eps}$, $t=t(x)=x^{\al +\al \eps -\eps \theta}$.
Let $F=\{  T^{(0)}_0 < t,  \tau^{(0)}_R > t \}$.
Then
\begin{align*} 
\bP^x(F^c) &\le \bP^x(  T^{(0)}_0 \ge  \tau^{(0)}_r)
 + \bP^x( \tau^{(0)}_r > t) +  \bP^x( \tau^{(0)}_R < t) \\
 &\le c (|x|/r)^\theta  + c r^\al /t + c t R^{-\al}.
\end{align*}
With the choices of $r$ and $t$ as above, one obtains
$\bP(F^c) \le 3 |x|^{-\eps \theta}$ provided $|x| \le R^{1/(1+\eps)}$.
Let $G=\{ Y_t \le x^{\al(1 +\eps)/2} \}$.
Then, using (\ref{ycheb}), we have
$\bP^x(G^c) \le  |x|^{-\eps \theta}$ also.

Suppose first that  $|x| \le R^{1/(1+\eps)}$. Then run
$X$ and $X^{(0)}$ until $S_1= T^{(0)}$.
We declare the run a success if both $F$ and $G$ occur,
so that success has a probability greater than
$1- c |x|^{-\eps \theta}$. If the run
is a success then we have $X_{S_1}=V_1$, where
$|V_1| \le  x^{\al(1 +\eps)/2}$.
We now repeat from the new starting point, and 
(if all the runs are successful) continue
until we obtain $X_{S_N}=V_N$ with $|V_N| \le R_1$.
Summing the probabilities of failures, we find that,
by choosing $R_1$ large enough, this can be made as small
as we like.

If we start at a point in $(-R_1, R_1)$,  a variant
of the argument above gives that, with probability $p_1>0$,
$X^{(0)}$ hits 0 before the first jump of $X^{(1)}$.
Finally, if  $R^{1/(1+\eps)}< |x| \le R/4$ then
running  $X^{(0)}$ until  $S_0= T^{(0)}$ we find with
probability $p_2>0$ that $|X_{S_0}|\le R^{1/(1+\eps)}$.

We deduce from this that $X$ satisfies $E_\al$ and EHI
with constants which do not depend on $R_1$. On the other
hand, $X$ only satisfies UJ$(\al)$ with a constant
of order $\log R_1$. This is enough to prove that the
`strong' form of the implication ``VD+ EHI + $E_\al
\Rightarrow$ PHI$(\al)$'' is false. That is (see Remark 
\ref{Rem-eff}), we cannot have  PHI$(\al)$ with a 
constant $C_P$ depending only on the constants in
VD, EHI and $E_\al$.

To actually obtain a single graph which satisfies
VD, EHI and $E_\al$ but not PHI$(\al)$, one needs
to modify the example above as follows. 
Take a rapidly increasing sequence $R_n$, define
$J_n$ analogously to $J_1$, and let $J=J_0 + \sum_{n \ge 1} J_n$.
This clearly fails to satisfy UJ$(\al)$, and so  
PHI$(\al)$ must also fail. 
However, arguments similar to the above show that
$E_\al$ and EHI still hold.

\begin{rem} {\rm
A recent paper  \cite{BS} gives necessary and sufficient
conditions for EHI to hold for $\al$-stable processes 
in $\bR^d$ with L\'evy measure of the form
$$ \nu(dx) = |x|^{-d-\al} f(x/|x|)dx , $$
where $f: S^{d-1} \to \bR_+$ is bounded and symmetric.
The condition in  \cite{BS} appears rather weaker than
UJS. If (as one may expect) the results of  \cite{BS}
hold also for processes on $\bZ^d$, this would 
give another class of examples when VD, EHI and $E_\al$ 
hold, but PHI$(\al)$ fails. 
} \end{rem}

\vskip 0.1 truein

\noi MTB: Department of Mathematics, University of British Columbia, Vancouver
V6T 1Z2, Canada

\noi RFB: Department of Mathematics, University of Connecticut, Storrs, 
CT 06269-3009, USA

\noi TK: Research Institute for Mathematical 
Sciences, Kyoto University, Kyoto 606-8502, Japan    

\end{document}